\documentclass[english, pdflatex]{siam/siamart220329}
\pdfoutput=1
\usepackage[T1]{fontenc}
\usepackage[latin9]{inputenc}
\usepackage{array}
\usepackage{bm}
\usepackage{multirow}
\usepackage{amsbsy}
\usepackage{amstext}
\usepackage{amssymb}
\usepackage{graphicx}
\usepackage{rotating}

\makeatletter

\providecommand{\tabularnewline}{\\}

\usepackage{amssymb}
\usepackage{amsmath}
\usepackage{bm}
\usepackage{graphicx}
\usepackage{caption}
\usepackage{subcaption}
\usepackage{tikz,pgfplots}

\pgfplotsset{compat=1.16}
\usetikzlibrary{external}
\usepackage{multirow, multicol}

\newcommand{\identity}{\mathbb{I}^{s \times s}}

\newcommand{\RR}{\mathbb{R}}
\newcommand{\NN}{\mathbb{N}}

\newcommand{\bs}[1]{\boldsymbol{#1}}

\makeatother

\usepackage{babel}

\begin{document}
\title{\markboth{M.R. Clines, V.E. Howle, and K.R. Long}{Efficient Preconditioners for IRK and IRKN Methods}Efficient
Order-Optimal Preconditioners for Implicit Runge--Kutta and Runge--Kutta--Nyström
Methods Applicable to a Large Class of Parabolic and Hyperbolic PDEs}
\author{Michael R. Clines\thanks{(Corresponding author) Department of Mathematics and Statistics, Texas
Tech University, Lubbock Texas, USA}, \and Victoria E. Howle\thanks{Department of Mathematics and Statistics, Texas Tech University, Lubbock
Texas, USA}, \and Katharine R. Long\thanks{Department of Mathematics and Statistics, Texas Tech University, Lubbock
Texas, USA}}
\maketitle
\begin{abstract}
We generalize previous work by Mardal\emph{, }Nilssen, and Staff (2007,
SIAM \emph{J. Sci. Comp. }v. 29\emph{, }pp. 361-375)\emph{ }and Rana,
Howle, Long, Meek, and Milestone (2021, SIAM \emph{J. Sci. Comp. }v.\emph{
}43\emph{, }pp. 475-495) on order-optimal preconditioners\emph{ }for
parabolic PDEs to a larger class of differential equations and methods.
The problems considered are those of the forms $u_{t}=-\mathcal{K}u+g$
and $u_{tt}=-\mathcal{{K}}u+g$, where the operator $\mathcal{{K}}$
is defined by $\mathcal{{K}}u:=-\nabla\cdot\left(\alpha\nabla u\right)+\beta u$
and the functions $\alpha$ and $\beta$ are restricted so that $\alpha>0$,
and $\beta\ge0$. The methods considered are A-stable implicit Runge--Kutta
methods for the parabolic equation and implicit Runge--Kutta--Nyström
methods for the hyperbolic equation. We prove the order optimality
of a class of block preconditioners for the stage equation system
arising from these problems, and furthermore we show that the LD and
DU preconditioners of Rana \emph{et al. }are in this class. We carry
out numerical experiments on several test problems in this class ---
the 2D diffusion equation, Pennes bioheat equation, the wave equation,
and the Klein--Gordon equation, with both constant and variable coefficients.
Our experiments show that these preconditioners, particularly the
LD preconditioner, are successful at reducing the condition number
of the systems as well as improving the convergence rate and solve
time for GMRES applied to the stage equations. 
\end{abstract}



\section{Introduction}

Implicit Runge--Kutta (IRK) methods are efficient time integration
methods that avoid the Dahlquist barriers that limit the order and
stability of implicit multistep methods. A disadvantage of IRK methods
has been the lack of effective preconditioners for the resulting linear
systems; however, some progress has been made. In \cite{mardal2007},
Mardal \emph{et al. }proved order optimality of a class of block diagonal
preconditioners for parabolic equations, where order optimality means
that the condition number of the preconditioned system is bounded
independently of the timestep $h_{t}$ and the mesh size $h_{x}$.
In \cite{masud2020}, Rana \emph{et al. }developed a preconditioner
based on a block LDU factorization that showed excellent performance
on the diffusion and advection-diffusion equations; however, order
optimality of their preconditioner was not proved. 

In this paper, we generalize the analysis of \cite{mardal2007} to
prove order optimality for a class of preconditioners for IRK systems
arising from a family of parabolic and hyperbolic partial differential
equations that includes the diffusion, Pennes bioheat\cite{pennes1948analysis},
wave, and Klein-Gordon equations. In particular, we prove that the
LD and DU preconditioners of Rana \emph{et al. }are order optimal
for all equations in this family. 

Let $\Omega$ be a bounded polygonal region in $\RR^{d}$ with boundary
$\partial\Omega$ partitioned into disjoint subsets $\partial\Omega_{N}$
and $\partial\Omega_{D}$. Let $\mu$ be either $1$ or $2$, and
let $\alpha\left(\mathbf{r}\right)$, $\beta\left(\mathbf{r}\right)$,
$g\left(\mathbf{r},t\right)$, $u_{0}\left(\mathbf{r}\right)$, and
$\dot{u}_{0}\left(\mathbf{r}\right)$ be real-valued functions, with
$\alpha$ and $\beta$ subject to the conditions $\alpha\left(\mathbf{r}\right)>0,\;\beta\left(\mathbf{r}\right)\ge0$
$\forall\mathbf{r}\in\Omega$. We define the linear differential operator
$\mathcal{{K}}$ by 
\[
\mathcal{{K}}u:=-\nabla\cdot\left(\alpha\nabla u\right)+\beta u.
\]
We focus on the linear initial-boundary value problem 
\begin{align}
\frac{{\partial^{\mu}u}}{\partial t^{\mu}} & =-\mathcal{{K}}u+g\;\;\;\text{in }\Omega\;\text{for }t>0\label{eq:PDE}\\
\alpha\boldsymbol{\hat{n}}\cdot\nabla u & =0\;\;\text{on \ensuremath{\partial\Omega_{N}}; \ensuremath{u=0}\;\;on }\Omega_{D}\label{eq:BCs}
\end{align}
with initial conditions
\begin{equation}
u\left(\mathbf{r},0\right)=u_{0}\left(\mathbf{r}\right)\label{eq:IC0}
\end{equation}
and, if $\mu=2$,
\begin{equation}
u_{t}\left(\mathbf{r},0\right)=\dot{u}_{0}\left(\mathbf{r}\right).\label{eq:IC1}
\end{equation}
With $\mu=1$ we have the diffusion ($\beta=0$) and Pennes bioheat
$\left(\beta>0\right)$ equations, while with $\mu=2$ we have the
wave $\left(\beta=0\right)$ and Klein--Gordon ($\beta>0$) equations. 

Time discretization will be done with an IRK method when $\mu=1$,
and with an implicit Runge--Kutta--Nyström (IRKN) method when $\mu=2$;
these methods are briefly reviewed in section \ref{subsec:Implicit-Runge-Kutta-methods}.
In subsection \ref{subsec:A-unified-formulation} we present a unified
formulation of the Galerkin weak forms of the stage equations for
IRK and IRKN methods. From this we develop a continuous boundary value
problem $\mathcal{{A}}\mathbf{k}=\mathbf{f}$ that must be solved
at every timestep for the stage variables $\mathbf{k}$; when this
problem is discretized in space with finite elements, we obtain a
linear system of equations. It is that linear system we are concerned
with preconditioning. 

Following \cite{mardal2007}, in section \ref{sec:Analysis-of-the}
we study the continuous form of the stage equations to show that the
mapping $\mathcal{{A}}$ is an isomorphism, and we apply that analysis
to study order-optimality of block preconditioners. Results of numerical
experiments are presented in section \ref{sec:NumericalResults},
and then conclusions and future directions are discussed in section
\ref{sec:Conclusions}.

\section{\label{sec:Time-and-Space}Time and space discretization}

For the remainder of this paper, vector-valued functions and spaces
of vector-valued functions are denoted by bold face symbols. The $L^{2}$
inner product and its induced norm are defined as usual as
\[
\left\langle u,v\right\rangle _{L^{2}}:=\int_{\Omega}uv\,d\Omega;\;\;\;\;\left\Vert u\right\Vert _{L^{2}}:=\sqrt{\left\langle u,u\right\rangle _{L_{2}}}.
\]
For a Sobolev inner product and norm we take
\[
\left\langle u,v\right\rangle _{H_{0}^{1}}:=\int_{\Omega}\alpha\nabla u\cdot\nabla v+\beta uv\,d\Omega;\;\;\;\;\left\Vert u\right\Vert _{H_{0}^{1}}:=\sqrt{\left\langle u,u\right\rangle _{H_{0}^{1}}},
\]
where $h_{t}>0$ is a timestep. In the corner case where $\beta=0$
and $\Omega_{D}=\emptyset$, our Sobolev inner product is not a true
inner product, and its norm is only a seminorm; however, we will soon
see that this is ultimately irrelevant because these operations won't
be used alone, so we tolerate this abuse of notation. With a number
of Runge--Kutta stages $s\in\NN$, we define the spaces
\begin{align*}
V & :=L^{2}(\Omega)\cap H_{0}^{1}(\Omega),\\
\mathbf{V} & :=V^{s}=(V\times\cdots\times V)=\mathbf{L}^{2}(\Omega)\cap h_{t}^{\mu}\mathbf{H}_{0}^{1}(\Omega),
\end{align*}
where $\otimes$ is the Cartesian product of vector spaces. The inner
products and norms in these spaces are
\begin{align*}
\left\langle u,v\right\rangle _{V} & :=\left\langle u,v\right\rangle _{L^{2}}+h_{t}^{\mu}\left\langle u,v\right\rangle _{H_{0}^{1}};\;\;\;\;\left\Vert u\right\Vert _{V}=\sqrt{\left\Vert u\right\Vert _{L^{2}}^{2}+h_{t}^{\mu}\left\Vert u\right\Vert _{H_{0}^{1}}^{2}}\\
\left\langle \mathbf{u},\mathbf{v}\right\rangle _{\mathbf{V}} & :=\left\langle \mathbf{u},\mathbf{v}\right\rangle _{\mathbf{L}^{2}}+h_{t}^{\mu}\left\langle \mathbf{u},\mathbf{v}\right\rangle _{\mathbf{H}_{0}^{1}};\;\;\;\;\left\Vert \mathbf{u}\right\Vert _{\mathbf{V}}=\sqrt{\left\Vert \mathbf{u}\right\Vert _{\mathbf{L}^{2}}^{2}+h_{t}^{\mu}\left\Vert \mathbf{u}\right\Vert _{\mathbf{H}_{0}^{1}}^{2}}.
\end{align*}
These inner products and norms on $V$ and $\mathbf{V}$ are true
inner products and norms even in the corner case where $\left\Vert \cdot\right\Vert _{H_{0}^{1}}$
is only a seminorm. Note that for the $L^{2},H_{0}^{1},$ and $V$
inner products, the vectorized inner products $\left\langle \mathbf{u},\mathbf{v}\right\rangle $
are related to the inner products $\left\langle u,v\right\rangle $
by
\[
\left\langle \mathbf{u},\mathbf{v}\right\rangle =\sum_{i=1}^{s}\left\langle u_{i},v_{i}\right\rangle .
\]
We will remove the space specification whenever it is clear which
inner product is being used. The operator $\mathcal{{K}}$ applied
to a vector $\mathbf{u}$ is to be interpreted such that $\left(\mathcal{K}\mathbf{u}\right)_{i}=\mathcal{{K}}u_{i}$.

\subsection{\label{subsec:Implicit-Runge-Kutta-methods}Implicit Runge--Kutta
methods}

An IRK method with $s$ stages is specified by its Butcher coefficients
$A\in\RR^{s\times s}$, $b\in\RR^{s}$, $c\in\RR^{s}$. Throughout
this paper we assume that we are working with an A-stable IRK or IRKN
method. For all such methods, the matrix $A$ will be nonsingular
and irreducible. With IRKN methods based on indirect collocation (see,
\emph{e.g.}, \cite{hairer1979,houwen1991}), the Butcher matrix $A$
is formed as $A=\hat{A}\hat{A}$, where $\hat{A}$ is the matrix from
a base IRK method. Our analysis will rely on the matrix $A$ being
weakly positive definite \cite{mardal2007,nilssen2005}. We have confirmed
this assumption through numerical experiments by examining the spectrum
of all the IRK and IRKN methods considered here.

Since the focus of timestepping methods is on the time variable, in
the following discussion we don't show explicitly the dependence on
spatial position $\mathbf{r}$; for example, we write $g\left(\mathbf{r},t\right)$
as $g\left(t\right)$ and $u^{n}\left(\mathbf{r}\right)$ simply as
$u^{n}$. 

The numerical solution at step $n-1$ is advanced to step $n$ by
the step formula
\begin{align*}
t^{n} & =t^{n-1}+h_{t}\\
u^{n} & =u^{n-1}+h_{t}\sum_{i=1}^{s}b_{i}k_{i}^{n},
\end{align*}
where the stage variables $k_{i}^{n}:\Omega\to\RR$, $i=1:s$, are
computed by solving the stage equations
\[
k_{i}^{n}=-\mathcal{{K}}\left(u^{n-1}+h_{t}\sum_{j=1}^{s}a_{ij}k_{j}^{n}\right)+g\left(t^{n-1}+c_{i}h_{t}\right);\;\;i=1:s,
\]
with homogeneous boundary conditions on each $k_{i}^{n}$. The stage
equations are a set of $s$ coupled boundary value problems. 

\subsubsection{\label{subsec:IRK-Nystr=0000F6m-methods}IRK--Nyström methods}

Runge--Kutta--Nyström methods are an extension of Runge--Kutta
methods suitable for second-order differential equations. A second
set of weights, $b'\in\RR^{s}$, is used so that the step formulas
for advancing the solution $u^{n}$ and the time derivative $\dot{u}^{n}$
are
\begin{align*}
{\displaystyle u^{n}} & =u^{n-1}+h_{t}\dot{u}^{n-1}+h_{t}^{2}\sum_{i=1}^{s}b_{i}k_{i}^{n}\\
\dot{u}^{n} & =\dot{u}^{n-1}+h_{t}\sum_{i=1}^{s}b_{i}^{\prime}k_{i}^{n}.
\end{align*}
The stage variables $k_{i}^{n}$ are found by solving 
\[
k_{i}^{n}=-\mathcal{{K}}\left(u^{n-1}+h_{t}c_{i}\dot{u}^{n-1}+h_{t}^{2}\sum_{j=1}^{s}a_{ij}k_{j}^{n}\right)+g\left(t^{n-1}+c_{i}h_{t}\right);\;\;i=1:s,
\]
with homogeneous boundary conditions. Note that although the equation
is second order in time, with a RKN method there are only $s$, not
$2s$, stage equations to be solved. 

\subsubsection{\label{subsec:A-unified-formulation}A unified formulation of the
IRK and IRKN stage equations}

Observe that we can easily express the stage equations for both IRK
(with $\mu=1$) and IRKN ($\mu=2$) in the unified form
\[
k_{i}^{n}=-\mathcal{{K}}\left(u^{n-1}+\left(\mu-1\right)h_{t}c_{i}\dot{u}^{n-1}+h_{t}^{\mu}\sum_{j=1}^{s}a_{ij}k_{j}^{n}\right)+g\left(t^{n-1}+c_{i}h_{t}\right).
\]
The time derivatives $\dot{u}^{n}$ are never needed when $\mu=1$.
If we introduce the variable
\begin{equation}
f_{i}^{n}=g\left(t^{n}+c_{i}h_{t}\right)+u^{n-1}+\left(\mu-1\right)h_{t}c_{i}\dot{u}^{n-1},\label{eq:RHS}
\end{equation}
then the $i$-th stage equation becomes
\begin{equation}
k_{i}^{n}+h_{t}^{\mu}\sum_{j=1}^{s}a_{ij}\mathcal{{K}}k_{j}^{n}=f_{i}^{n}.\label{eq:stage-eqns}
\end{equation}
As a final notational simplification, introduce the vector of stage
variables 
\[
\mathbf{k}^{n}=\left[\begin{array}{cccc}
k_{1}^{n} & k_{2}^{n} & \cdots & k_{s}^{n}\end{array}\right]^{T},
\]
the vector $\mathbf{f}^{n}=\left[\begin{array}{cccc}
f_{1}^{n} & f_{2}^{n} & \cdots & f_{s}^{n}\end{array}\right]^{T}$, the $s\times s$ identity matrix $\identity$, and the identity
operator $\mathcal{{I}}:V\to V$. Then the system of stage equations
is written compactly as
\[
\left(\identity\otimes\mathcal{{I}}+h_{t}^{\mu}A\otimes\mathcal{{K}}\right)\mathbf{k}^{n}=\mathbf{f}^{n},
\]
where $\otimes$ is the Kronecker product. We can also define $\mathcal{{A}}:\mathbf{V}\to\mathbf{V}^{*}$
as
\begin{equation}
\mathcal{{A}}:=\identity\otimes\mathcal{{I}}+h_{t}^{\mu}A\otimes\mathcal{{K}},\label{eq:cont-A}
\end{equation}
in which case the equation to be solved is simply 
\begin{equation}
\mathcal{{A}}\mathbf{k}^{n}=\mathbf{f}^{n}.\label{eq:cont-A-eqn}
\end{equation}

\subsection{Finite element discretization}

We now develop a variational form of the stage equations. Begin by
observing that for all $v\in V$, we have
\[
\int_{\Omega}v\mathcal{{K}}u\,d\Omega=\int_{\Omega}\alpha\nabla u\cdot\nabla v+\beta uv\,d\Omega=\left\langle v,u\right\rangle _{H_{0}^{1}}.
\]
Now assert an orthogonal residual condition with test function $v_{i}\in V$
on equation \ref{eq:stage-eqns}, obtaining
\begin{equation}
\left\langle k_{i}^{n},v_{i}\right\rangle _{L^{2}}+h_{t}^{\mu}\sum_{j=1}^{s}a_{ij}\left\langle k_{j}^{n},v_{i}\right\rangle _{H_{0}^{1}}=\left\langle f_{i}^{n},v_{i}\right\rangle _{L^{2}}.\label{eq:expanded-weak}
\end{equation}

Choose an $N$-dimensional approximating subspace $\mathbf{V}_{h}\subset\mathbf{V}$,
where $h>0$ is a notational representation of the mesh size. Quantities
with subscript $h$ are discretized in this space; for example, the
function $k_{i,h}^{n}\in V_{h}$ will be the $i$-th discrete stage
variable at timestep $n$. The vectorized discrete space is $\mathbf{V}_{h}=\left(V_{h}\right)^{s}$.
Let $\left\{ \phi_{\ell}\right\} _{\ell=1}^{N}$ be a basis for $V_{h}$,
and form the matrices $M\in\RR^{N\times N}$ and $F\in\RR^{N\times N}$
with entries
\begin{align}
M_{\ell m} & :=\left\langle \phi_{\ell},\phi_{m}\right\rangle _{L^{2}}\label{eq:mass-mtx}\\
F_{\ell m} & :=\left\langle \mathcal{{K}}\phi_{m},\phi_{\ell}\right\rangle _{L^{2}}=\left\langle \alpha\nabla\phi_{m},\nabla\phi_{\ell}\right\rangle _{L^{2}}+\left\langle \beta\phi_{m},\phi_{\ell}\right\rangle _{L^{2}}=\left\langle \phi_{m},\phi_{\ell}\right\rangle _{H_{0}^{1}}\label{eq:F-mtx-2}
\end{align}
and the vector $f_{i,h}^{n}\in\RR^{N}$ with entries
\begin{equation}
\left(f_{i,h}^{n}\right)_{\ell}=\left\langle f_{i}^{n},\phi_{\ell}\right\rangle _{L^{2}}.\label{eq:f-vec}
\end{equation}
The $i$-th stage equation from (\ref{eq:expanded-weak}) then becomes
\[
Mk_{i,h}^{n}+h_{t}^{\mu}\sum_{j=1}^{n}a_{ij}Fk_{j,h}^{n}=f_{i,h}^{n},
\]
and the full system in vectorized notation is
\begin{equation}
\left(\identity\otimes M+h_{t}^{\mu}A\otimes F\right)\mathbf{k}_{h}^{n}=\mathbf{f}_{h}^{n}.\label{eq:discr-sys}
\end{equation}
The matrix on the left hand side is the discretization of the operator
$\mathcal{{A}}$, so we denote it as $\mathcal{{A}}_{h}$,
\[
\mathcal{{A}}_{h}=\identity\otimes M+h_{t}^{\mu}A\otimes F.
\]

\section{\label{sec:Analysis-of-the}Analysis of the continuous stage equations
and preconditioner}

We now analyze the continuous form of the stage equations (\ref{eq:cont-A-eqn})
\[
\mathcal{A}\mathbf{k}^{n}=\mathbf{f}^{n},
\]
where the operator $\mathcal{A}:\mathbf{V}\mapsto\mathbf{V}^{*}$
has been defined in (\ref{eq:cont-A}),
\begin{align*}
\mathcal{A} & :=\identity\otimes\mathcal{{I}}+h_{t}^{\mu}A\otimes\mathcal{{K}}.
\end{align*}
We note the similarity to the stage equation operator $\mathcal{A}:=I-h_{t}A\otimes\nabla^{2}$
considered in \cite{mardal2007}. We therefore follow the path laid
out by \cite{mardal2007}: we establish that $\mathcal{{A}}$ is an
isomorphism given certain assumptions about the Butcher coefficient
matrix $A$, and obtain the following theorem.
\begin{theorem}
\label{thm:iso-A}Let $A\in\RR^{s\times s}$ be weakly positive definite,
and let $h_{t}>0$, $\mu\in\text{\ensuremath{\left\{  1,2\right\} } }$.
Then the operator $\mathcal{A}=\identity\otimes\mathcal{{I}}+h_{t}^{\mu}A\otimes\mathcal{{K}}$
is an isomorphism, and there exist finite constants $c_{1}$,$c_{2}$,
independent of the timestep $h_{t}$, such that
\begin{align*}
\left\Vert \mathcal{{A}}\right\Vert _{\mathcal{{L\left({\bf V,}{\bf V^{*}}\right)}}} & \leq c_{1}\\
\left\Vert \mathcal{{A}}^{-1}\right\Vert _{\mathcal{\mathcal{{L}\left({\bf V^{*},{\bf V}}\right)}}} & \leq c_{2}.
\end{align*}
\end{theorem}

Because the proof is a straightforward extension of that found in
\cite{mardal2007}, we leave the proof to Appendix \ref{sec:Proof-of-theorem}. 

Now consider a continuous preconditioner of the form 
\begin{equation}
\mathcal{P}=I+h_{t}^{\mu}P\otimes\mathcal{K},\label{eq:cont-prec}
\end{equation}
where $P$ is a weakly positive definite matrix; this is simply the
operator $\mathcal{{A}}$ but with the Butcher matrix $A$ replaced
by an approximation $P$ that has preserved the weakly positive definiteness
of $A$. From Theorem \ref{thm:iso-A} we see immediately that $\mathcal{{P}}$
is an isomorphism. Since $\mathcal{P}$ is an isomorphism between
$\mathbf{V}$ and $\mathbf{V}^{*}$, its inverse $\mathcal{P}^{-1}$
exists and maps $\mathbf{V}^{*}$ to $\mathbf{V}$, and from Theorem
\ref{thm:iso-A} we have the bounds
\[
\left\Vert \mathcal{{P}}\right\Vert \leq d_{1};\;\;\;\;\left\Vert \mathcal{{P}}^{-1}\right\Vert \leq d_{2}.
\]
Furthermore, the composition $\mathcal{{P}}^{-1}\mathcal{{A}}:\mathbf{V}\to\mathbf{V}$
is also an isomorphism with bounds 
\[
\left\Vert \mathcal{P}^{-1}\mathcal{A}\right\Vert \leq c_{1}d_{2}\mbox{ and }\left\Vert \mathcal{A}^{-1}\mathcal{P}\right\Vert \leq c_{2}d_{1}.
\]
From this we immediately establish a bound on the condition number
of the continuous left-preconditioned operator $\mathcal{P}^{-1}\mathcal{A}$
that is independent of our choice of $h_{t}$, that is,
\[
\kappa(\mathcal{P}^{-1}\mathcal{A})=\left\Vert \mathcal{{P}}^{-1}\mathcal{{A}}\right\Vert \left\Vert \mathcal{{A}}^{-1}\mathcal{{P}}\right\Vert \leq c_{1}d_{2}c_{2}d_{1}.
\]
Hence this family of preconditioners is a set of order optimal preconditioners
with respect to time. Next we establish that these preconditioners
are also order optimal with respect to the discretization parameter
$h$. 

We now consider the following discrete subspaces 
\[
V_{h}\subset V\quad\mbox{and}\quad\mathbf{V}_{h}=(V_{h})^{s}\subset\mathbf{V}.
\]
The discretized versions $\mathcal{A}_{h}$ and $\mathcal{{P}}_{h}$
of the above continuous operators $\mathcal{A}$ and $\mathcal{P}$
are given by 
\begin{align}
\mathcal{A}_{h} & :=\identity\otimes M+h_{t}^{p}A\otimes F\label{eq:A-disc}\\
\mathcal{P}_{h} & :=\identity\otimes M+h_{t}^{p}P\otimes F,\label{eq:P-disc}
\end{align}
where $F,M$ are the stiffness and mass matrices, respectively. Since
$\mathbf{V}_{h}\subset\mathbf{V}$, we obtain 
\[
\kappa(\mathcal{P}_{h}^{-1}\mathcal{A}_{h})\leq\kappa(\mathcal{P}^{-1}\mathcal{A})\leq c_{1}d_{2}c_{2}d_{1},
\]
showing that $\mathcal{P}^{-1}\mathcal{A}$ is bounded independently
of the discretization parameter $h$. We conclude that the preconditioner
$\mathcal{{P}}$ is order optimal with respect to both timestep $h_{t}$
and spatial discretization parameter $h$. 

Now that we have established that any preconditioner $\mathcal{{P}}$
with $P$ weakly positive definite is order-optimal, we now focus
on several examples of interest. We can define a few particular instances
of this preconditioner based on our selection of $P$, specifically
where $P$ is a preconditioner for the Butcher coefficient matrix
$A$ arising from our choice of timestepper. 

The motivation behind this decision comes from the observation that
$\kappa(\mathcal{P}^{-1}\mathcal{A})\propto\kappa(P^{-1}A)$ for $\beta=0$,
as established in \cite{staff2006}, \cite{mardal2007}, and \cite{masud2020}.
In other words, if we can construct $\mathcal{P}$ such that $P$
is a good preconditioner for $A$, then $\mathcal{P}$ is likely to
be a good preconditioner for $\mathcal{A}$. Some simple examples
are the block Jacobi and block Gauss-Seidel preconditioners studied
by \cite{staff2006,mardal2007} and the LDU-factorization preconditioners
introduced in \cite{masud2020}. Additionally, we will also show the
results of using the upper triangular part of $A,$ namely $\mathcal{{P}}_{TRIU}$.

The order optimality of $\mathcal{{P}}_{LD}$ and $\mathcal{{P}}_{DU}$
was suggested experimentally, but not proved, by \cite{masud2020}.
All five of these preconditioners are easily constructed and applied;
the factorizations in the LD and DU preconditioners are done on the
small Butcher coefficient matrices $A$, not on the full matrix $\mathcal{{A}}_{h}$.
We have numerically verified the weak positive definiteness of all
the above preconditioners. See Appendix \ref{sec:Preconditioner-Def}
for the specific definition of each preconditioner.

\section{\label{sec:NumericalResults} Numerical Experiments}

We conduct a series of numerical experiments to investigate the behavior
of several preconditioners in problems of the type considered in this
paper. Left preconditioning is used throughout. In all experiments,
homogeneous Neumann boundary conditions are used and spatial discretization
is done with Galerkin finite elements using first-degree Lagrange
basis functions on a 2D triangular mesh. For each problem, we consider
both constant and non-constant coefficients. For the constant coefficient
experiments, we set $\alpha=\beta=1$, and for the variable coefficient
experiments we use

\[
\begin{aligned}\alpha\left(\mathbf{r}\right):= & 1+0.2xy\\
\beta\left(\mathbf{r}\right):= & 1+0.3\sin(\pi x)\cos(\pi y)
\end{aligned}
.
\]
On the domain $\Omega=[-1,1]\times[-1,1]$ these functions are strictly
positive. 

For parabolic problems, the Radau IIA method is used; for the hyperbolic
problems, an IRKN method based on the Gauss--Legendre method is used.
Unless otherwise specified, the timestep $h_{t}$ is chosen to depend
on the mesh size $h$ so that the spatial interpolation error and
temporal global truncation error are of comparable magnitude \cite{masud2020};
we set $h_{t}=$ $h^{\frac{{p+1}}{q}}$ where $p$ is the degree of
the basis polynomials, $q$ is $2s$ for Gauss--Legendre timestepping
and $q=2s-1$ for Radau IIA timestepping. 

In GMRES calculations, the method of manufactured solutions is used
to construct exact solutions so that we can calculate relative errors
as well as relative residuals. 

The matrices $M$ and $F$ are assembled using the Sundance finite
element toolkit \cite{long2010}; when doing solves, the full matrices
$\mathcal{A}_{h}$ are never assembled as all matrix-vector multiplications
and diagonal subsolves can be carried out using only $M$ and $F$.
Other than the assembly of $M$ and $F$, all computations were implemented
in Python using the open source packages SciPy\cite{scipy} and PyAMG\cite{pyamg}
on a machine with an AMD Ryzen 7 3800 3.89GHz processor and 16GB of
RAM. 

\subsection{Constant coefficient problems\label{subsec:Constant-Coefficients-1}}

Our first set of calculations is done using constant coefficients.
We compute condition numbers for the constant coefficient preconditioned
wave and Klein--Gordon equations as functions of $h$, $h_{t}$ and
$s$ and for a variety of preconditioners (subsection \ref{subsec:Wave-Equation});
eigenvalue spectra for the preconditioned constant-coefficient wave
equation (subsection \ref{subsec:Spectrum-of-preconditioned}); and
field of values for the preconditioned constant-coefficient Klein--Gordon
equation (subsection \ref{subsec:Field-of-values}). 

Condition numbers for the preconditioned diffusion equation with constant
coefficients have already been investigated extensively in \cite{mardal2007,staff2006,masud2020},
so we don't repeat that analysis here. We defer experiments on the
diffusion and Pennes equations until our investigation of variable
coefficient problems. 

\subsubsection{Conditioning for the wave and Klein-Gordon equations\label{subsec:Wave-Equation}}

Table \ref{table:left-cn-wave-1} shows calculated condition numbers
of the system matrix $\mathcal{{A}}_{h}$ arising from the discretized
wave equation, constructed with mesh sizes $h=2^{-4},\:2^{-5}$ and
then preconditioned with $\mathcal{{P}}_{J},\,\mathcal{{P}}_{GSL}$,~$\mathcal{\mathcal{{P}}}_{LD}$,~$\mathcal{\mathcal{{P}}}_{DU}$,
and $\mathcal{\mathcal{{P}}}_{TRIU}$. The time step is defined to
be $h_{t}=h^{\frac{2}{2s}}$, as specified above. The timestepping
methods chosen are the s-stage IRKN Gauss--Legendre methods where
$s$ varies from 2--5.

For $s=2$, all preconditioners reduce the condition number by a factor
between 2.5 and 50. As $s$ is increased, we see a sharp rise in the
condition number of the original system, and we see an increase in
the condition number of the systems preconditioned with $\mathcal{\mathcal{{P}}}_{J},$$\mathcal{\mathcal{{P}}}_{GSL}$,
$\mathcal{P}_{DU}$, and $\mathcal{\mathcal{{P}}}_{TRIU}$. However,
$P_{LD}$ maintains a consistently small condition number, growing
only by $\sim4-5$ as stage number increases and staying under 10
for all $h$ and $s$ considered.
\begin{table}
\centering{}%
\begin{tabular}{l|r@{\extracolsep{0pt}.}lr@{\extracolsep{0pt}.}lr@{\extracolsep{0pt}.}lr@{\extracolsep{0pt}.}lr@{\extracolsep{0pt}.}lr@{\extracolsep{0pt}.}l}
\hline 
GL-(s) & \multicolumn{2}{c}{$\kappa(\mathcal{A})$} & \multicolumn{2}{c}{$\kappa(\mathcal{P}_{J}^{-1}\mathcal{A})$} & \multicolumn{2}{c}{$\kappa(\mathcal{P}_{GSL}^{-1}\mathcal{A})$} & \multicolumn{2}{c}{$\kappa(\mathcal{P}_{LD}^{-1}\mathcal{A})$} & \multicolumn{2}{c}{$\kappa(\mathcal{P}_{DU}^{-1}\mathcal{A})$} & \multicolumn{2}{c}{$\kappa(\mathcal{P}_{TRIU}^{-1}\mathcal{A})$}\tabularnewline
\hline 
\multicolumn{13}{c}{$h=2^{-4}$}\tabularnewline
\hline 
GL-2 & 68&22 & 17&43 & 4&39 & 2&06 & 7&68 & 20&46\tabularnewline
GL-3 & 410&88 & 389&36 & 43&73 & 4&80 & 233&82 & 595&84\tabularnewline
GL-4 & 966&05 & 395&42 & 656&31 & 5&41 & 314&39 & 481&09\tabularnewline
GL-5 & 1807&02 & 1023&39 & 753&04 & 8&01 & 1145&59 & 1377&70\tabularnewline
\hline 
\multicolumn{13}{c}{$h=2^{-5}$}\tabularnewline
\hline 
GL-2 & 105&40 & 18&91 & 4&55 & 2&12 & 8&15 & 22&07\tabularnewline
GL-3 & 802&25 & 413&53 & 45&88 & 6&03 & 244&96 & 630&85\tabularnewline
GL-4 & 2276&20 & 484&52 & 738&48 & 6&32 & 414&06 & 588&11\tabularnewline
GL-5 & 4676&51 & 1295&65 & 893&63 & 9&62 & 1432&78 & 1443&22\tabularnewline
\hline 
\end{tabular}\caption{Condition numbers of left-preconditioned matrices with preconditioners
Jacobi, Gauss--Seidel, $LD$, $DU$, and Upper Triangular, with spatial
discretization sizes $h=2^{-4}$,~$2^{-5}$ applied to a 2D wave
equation with constant coefficients with $s$-stage Gauss--Legendre.
The time-steps are given by $h_{t}$ = $h^{\frac{p+1}{q}}$ where
$p=1$ is the degree of the Lagrange polynomial basis functions in
space and $q$ is the order of the method (2$s$ for Gauss--Legendre).}
\label{table:left-cn-wave-1}
\end{table}
\begin{table}
\begin{centering}
\begin{tabular}{c|r@{\extracolsep{0pt}.}lr@{\extracolsep{0pt}.}lr@{\extracolsep{0pt}.}lr@{\extracolsep{0pt}.}lr@{\extracolsep{0pt}.}lr@{\extracolsep{0pt}.}l}
\hline 
GL-(s) & \multicolumn{2}{c}{$\kappa(\mathcal{A})$} & \multicolumn{2}{c}{$\kappa(\mathcal{P}_{J}^{-1}\mathcal{A})$} & \multicolumn{2}{c}{$\kappa(\mathcal{P}_{GSL}^{-1}\mathcal{A})$} & \multicolumn{2}{c}{$\kappa(\mathcal{P}_{LD}^{-1}\mathcal{A})$} & \multicolumn{2}{c}{$\kappa(\mathcal{P}_{DU}^{-1}\mathcal{A})$} & \multicolumn{2}{c}{$\kappa(\mathcal{P}_{TRIU}^{-1}\mathcal{A})$}\tabularnewline
\hline 
\multicolumn{13}{c}{$h_{t}=5.0$}\tabularnewline
\hline 
GL-2 & \multicolumn{2}{c}{24,776} & 29&03 & 8&31 & 3&38 & 9&65 & 36&08\tabularnewline
GL-3 & \multicolumn{2}{c}{70,447} & 1790&89 & 180&33 & 6&10 & 873&32 & 2857&00\tabularnewline
GL-4 & \multicolumn{2}{c}{251,624} & 3025&13 & 6495&80 & 33&89 & 1634&59 & 3621&96\tabularnewline
GL-5 & \multicolumn{2}{c}{469,410} & 12,133&55 & 63,830&62 & 82&16 & 10,919&53 & 15,574&36\tabularnewline
\hline 
\multicolumn{13}{c}{$h_{t}=.5$}\tabularnewline
\hline 
GL-2 & 3859&26 & 73&96 & 9&81 & 4&27 & 32&63 & 78&77\tabularnewline
GL-3 & 6267&10 & 786&89 & 89&95 & 11&77 & 534&50 & 1168&69\tabularnewline
GL-4 & 7644&23 & 653&50 & 926&44 & 9&11 & 75&06 & 796&45\tabularnewline
GL-5 & 9347&42 & 1524&38 & 9678&74 & 11&33 & 1915&84 & 1936&84\tabularnewline
\hline 
\multicolumn{13}{c}{$h_{t}=.005$}\tabularnewline
\hline 
GL-2 & 438&43 & 2&43 & 1&10 & 1&39 & 2&39 & 2&57\tabularnewline
GL-3 & 468&61 & 3&74 & 1&77 & 1&92 & 3&75 & 3&96\tabularnewline
GL-4 & 526&19 & 2&92 & 1&02 & 1&65 & 2&98 & 2&92\tabularnewline
GL-5 & 529&67 & 2&97 & 1&02 & 1&71 & 2&99 & 2&97\tabularnewline
\hline 
\end{tabular}
\par\end{centering}
\caption{Condition numbers of left-preconditioned matrices with preconditioners
Jacobi, Gauss--Seidel, $LD$, $DU$, and Upper Triangular, with time
step sizes $h_{t}=5.0,0.5,0.005$ and $h=2^{-5}$ applied to a 2D
Klein--Gordon equation with constant coefficients with $s$-stage
Gauss--Legendre.}

\label{table:left-cn-kg_time-1}
\end{table}

Table \ref{table:left-cn-kg_time-1} shows the results of preconditioning
the Klein--Gordon equation (with constant coefficients) with $h=2^{-5}$
held fixed, and with timesteps $h_{t}=5.0,\,0.5,\,0.005$ advanced
with Gauss--Legendre 2--5. All preconditioners reduce $\kappa$
by an order of magnitude or more, with $\mathcal{{P}}_{LD}$ superior
except at the smallest timesteps, where $\mathcal{{P}}_{GSL}$ is
slightly better. It is worth noting that as $h_{t}\to0,$ the original
system $\mathcal{{A}}$ converges to a block diagonal mass matrix,
which is easily handled by any preconditioner; indeed, all five preconditioners
considered perform well at $h_{t}=0.005$.

\subsubsection{\label{subsec:Spectrum-of-preconditioned}Spectrum of preconditioned
wave equation}

It is well known that when it comes to preconditioning the iterative
solver GMRES, the condition number alone does not necessarily predict
a preconditioner's effectiveness. 

\begin{figure}
\begin{centering}
\includegraphics[scale=0.75]{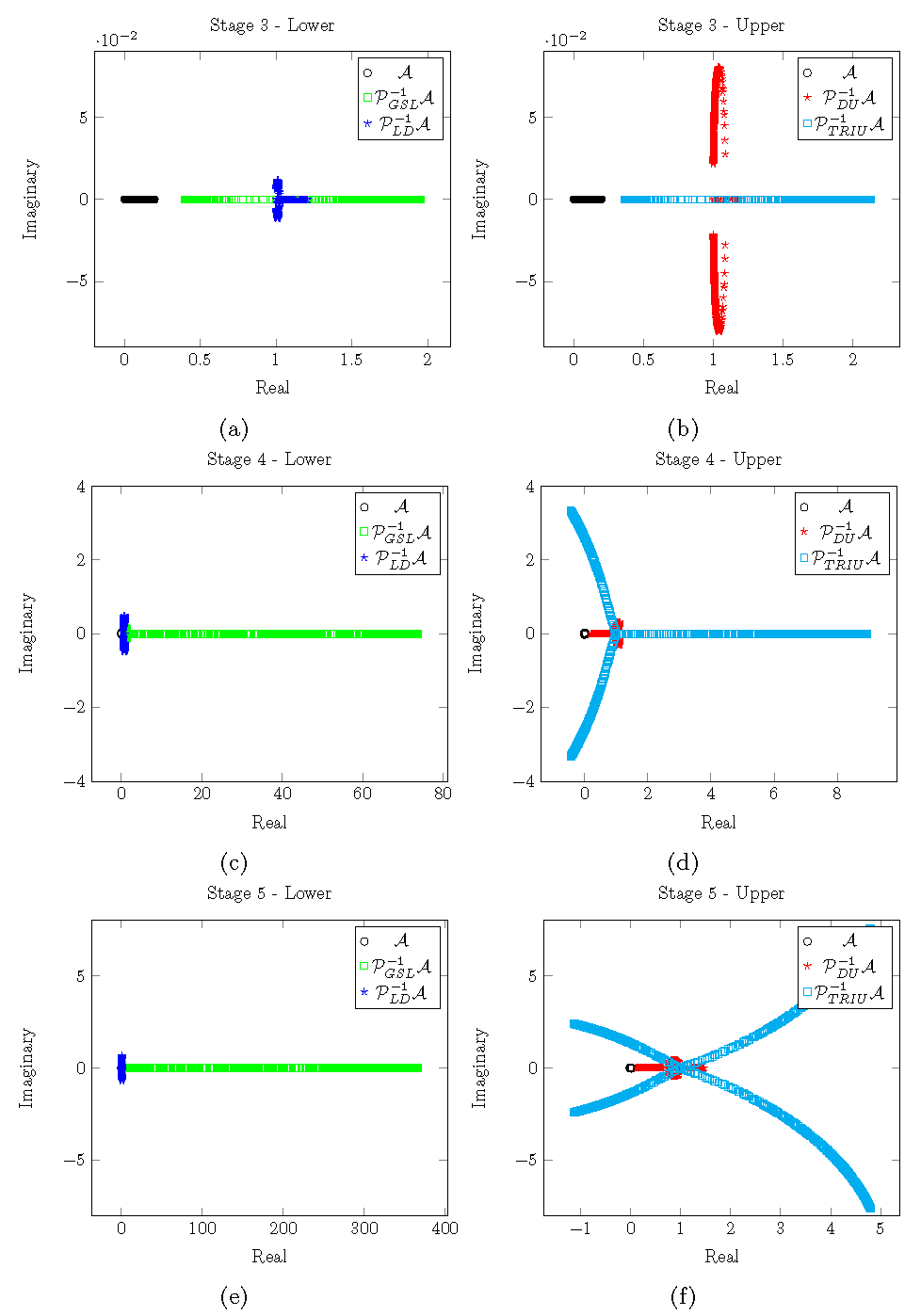}
\par\end{centering}
\caption{Eigenvalues of the left-preconditioned matrices arising from a 2D
wave equation, with $h=2^{-4}$ and using Gauss--Legendre with stages:
(a-b) $s=3$, (c-d) $s=4$ , and (e-f) $s=5$. The $x$-axis and $y$-axis
are the real and imaginary axes, respectively. The left column of
figures show eigenvalues for the lower triangular preconditioned systems,
and the right column of figures shows eigenvalues for the upper triangular
preconditioned systems.}

\label{Figure: SpectrumPlot-1-1}
\end{figure}

Generally, it is desirable to cluster the eigenvalues away from 0.
In Figure \ref{Figure: SpectrumPlot-1-1}, we show the spectrum of
the preconditioned wave equation with constant coefficients, spatial
discretization size $h=2^{-4}$, and time step $h_{t}=h^{\frac{{1}}{s}}$,
coupled with the Gauss--Legendre method of stages 3--5. The three
rows contain results for the three stages, with $s$ decreasing downwards.
The left column of figures shows the spectrum of the block lower triangular
preconditioned ($\mathcal{{P}}_{GSL}$ and $\mathcal{{P}}_{LD}$)
systems and the right shows the spectrum of the block upper triangular
($\mathcal{{P}}_{DU}$ and $\mathcal{{P}}_{TRIU}$) preconditioned
systems. The spectrum of the unpreconditioned system is also shown,
in black. As we can see in all the cases, the eigenvalues of the original
systems are clustered near 0 while the preconditioners $\mathcal{{P}}_{LD}$
and $\mathcal{{P}}_{DU}$ tend to cluster their eigenvalues near 1,
with LD achieving the better clustering of the two. 

\subsubsection{\label{subsec:Field-of-values}Field of values for preconditioned
Klein--Gordon equation}

\begin{figure}
\begin{centering}
\includegraphics[clip,scale=0.85]{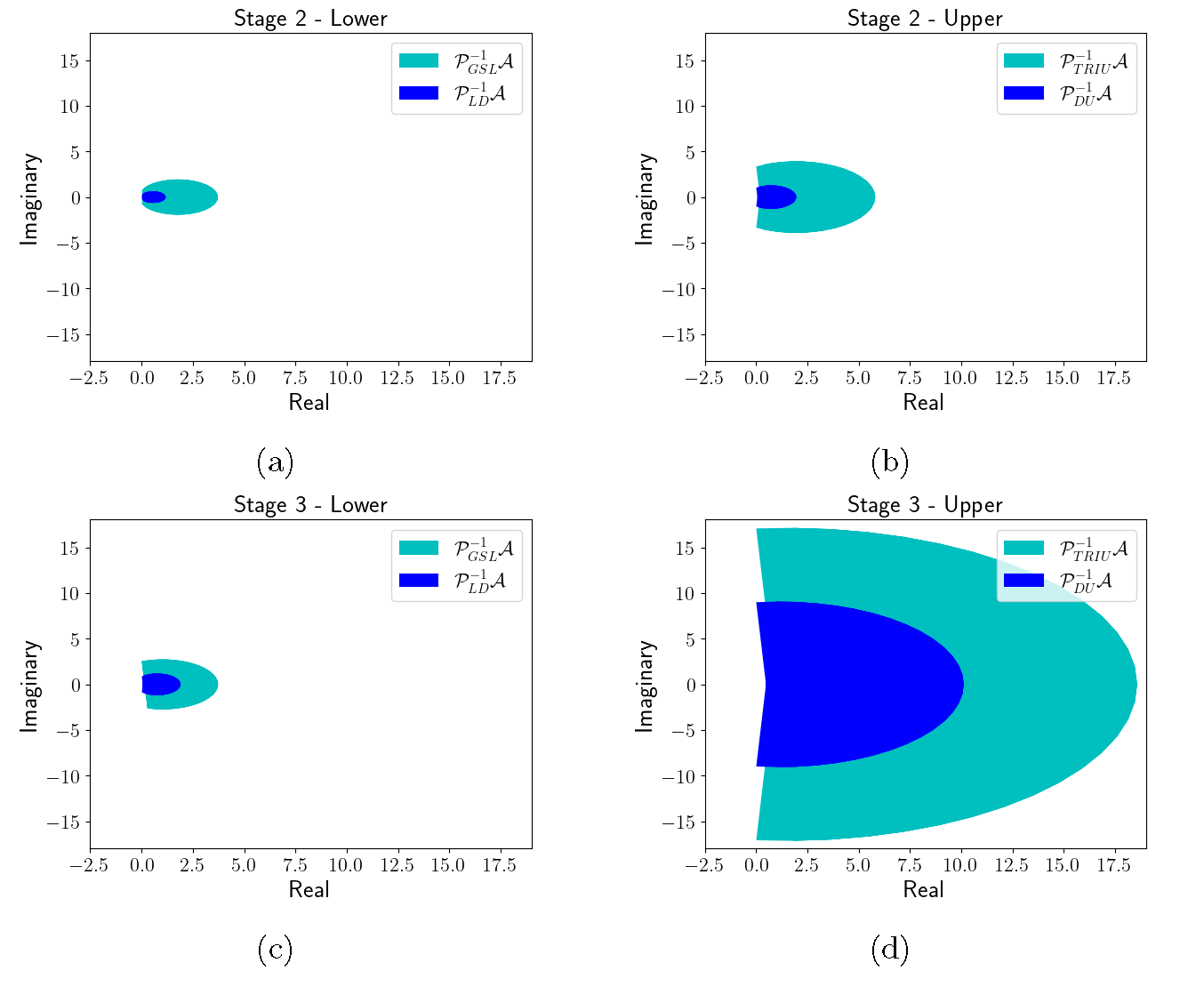}
\par\end{centering}
\caption{\label{fig:Numerical-Range-1}Numerical range in the complex plane
of preconditioned matrices arising from a Klein--Gordon equation
with simplified coefficients, with $h=2^{-4}$, and using Gauss--Legendre
with stages: (a-b) $s=2$ and (c-d) $s=3$. The left column of figures
shows the lower triangular preconditioned systems, and the right column
of figures shows the upper triangular preconditioned systems.}
\end{figure}

The field of values (FOV), or numerical range, can also indicate the
behavior of a GMRES preconditioner \cite{liesen2012}. Since one of
the standard worst-case error bounds on GMRES is given by the distance
of the boundary of the FOV to the origin \cite{liesen2012}, in Figure
\ref{fig:Numerical-Range-1} we show several plots depicting the numerical
range of $\mathcal{{A}}$ for the Klein--Gordon equation with constant
coefficients. The mesh size $h=2^{-4}$ and time step $h_{t}=h^{\frac{1}{s}}$
with Gauss--Legendre-2 and Gauss--Legendre-3 are used to create
the figures. It is clear from these results that $\mathcal{{P}}_{GSL}$
and $\mathcal{{P}}_{LD}$ have tighter numerical ranges (smaller distance
from zero) than the upper-triangular preconditioned systems, and that
$\mathcal{{P}}_{LD}$ produces a tighter numerical range than $\mathcal{P}_{GSL}$.

\subsection{\label{subsec:Variable-coefficient-problems}Variable coefficient
problems}

\subsubsection{\label{subsec:parab-cond}Condition numbers for the diffusion and
Pennes equations}

Next we consider conditioning of the preconditioned diffusion and
Pennes bioheat equations. In this set of experiments, we discretize
our domain with mesh size $h=2^{-4}$, and use time steps $h_{t}=h^{\frac{2}{2s-1}}$,
where $s$ is the stage number for the IRK Radau IIA method. As these
problems are parabolic, an L-stable method should be used. We follow
\cite{mardal2007,masud2020,staff2006} in choosing the Radau IIA method. 

Condition numbers were computed for these problems with the $GSL$,
$LD$, and $DU$ preconditioners. All three perform remarkably well,
with $LD$ performing consistently best, and with $GSL$ as a runner-up.
We also note that the variable coefficients appear to have a minimal
impact on the condition number, hence the small variations between
the two sets of data. 

\begin{table}
\centering{}%
\begin{tabular}{l|r@{\extracolsep{0pt}.}lr@{\extracolsep{0pt}.}lr@{\extracolsep{0pt}.}lr@{\extracolsep{0pt}.}l}
\hline 
RIIA & \multicolumn{2}{c}{$\kappa(\mathcal{A})$} & \multicolumn{2}{c}{$\kappa(\mathcal{P}_{GSL}^{-1}\mathcal{A})$} & \multicolumn{2}{c}{$\kappa(\mathcal{P}_{LD}^{-1}\mathcal{A})$} & \multicolumn{2}{c}{$\kappa(\mathcal{P}_{DU}^{-1}\mathcal{A})$}\tabularnewline
\hline 
\multicolumn{9}{c}{Diffusion ($\beta=0$)}\tabularnewline
\hline 
RIIA - 2 & 973&90 & 2&99 & 1&67 & 24&20\tabularnewline
RIIA - 3 & 3372&07 & 7&51 & 2&48 & 84&65\tabularnewline
RIIA - 4 & 5743&65 & 15&89 & 3&11 & 174&67\tabularnewline
RIIA - 5 & 8280&81 & 30&27 & 3&86 & 306&03\tabularnewline
\hline 
\multicolumn{9}{c}{Pennes ($\beta>0$)}\tabularnewline
\hline 
RIIA - 2 & 964&09 & 2&98 & 1&65 & 23&76\tabularnewline
RIIA - 3 & 3286&59 & 7&44 & 2&43 & 82&01\tabularnewline
RIIA - 4 & 5693&23 & 15&68 & 2&27 & 81&32\tabularnewline
RIIA - 5 & 8070&54 & 29&95 & 3&77 & 297&37\tabularnewline
\hline 
\end{tabular}\caption{Condition numbers of left-preconditioned matrices with preconditioners
Gauss--Seidel, $LD$, and $DU$ with spatial discretization size
$h=2^{-4}$ applied to 2D diffusion and Pennes equations with variable
coefficients with $s$-stage Radau IIA. The time-steps are given by
$h_{t}$ = $h^{\frac{2}{2s-1}}$.}
\label{table:left-parabolic}
\end{table}

\subsubsection{\label{subsec:hyp-cond}Condition numbers for the wave and Klein--Gordon
equations}

\begin{table}
\centering{}%
\begin{tabular}{c|r@{\extracolsep{0pt}.}lr@{\extracolsep{0pt}.}lr@{\extracolsep{0pt}.}lr@{\extracolsep{0pt}.}l}
\hline 
GL & \multicolumn{2}{c}{$\kappa(\mathcal{A})$} & \multicolumn{2}{c}{$\kappa(\mathcal{P}_{GSL}^{-1}\mathcal{A})$} & \multicolumn{2}{c}{$\kappa(\mathcal{P}_{LD}^{-1}\mathcal{A})$} & \multicolumn{2}{c}{$\kappa(\mathcal{P}_{DU}^{-1}\mathcal{A})$}\tabularnewline
\hline 
\multicolumn{9}{c}{$\beta=0$}\tabularnewline
\hline 
GL-2 & 167&67 & 5&50 & 2&72 & 15&04\tabularnewline
GL-3 & 691&58 & 74&71 & 8&32 & 467&92\tabularnewline
GL-4 & 1906&88 & 1049&11 & 9&70 & 994&91\tabularnewline
GL-5 & 3923&20 & 1776&26 & 16&09 & 3503&47\tabularnewline
\hline 
\multicolumn{9}{c}{$\beta>0$}\tabularnewline
\hline 
GL-2 & 167&61 & 5&50 & 2&71 & 15&01\tabularnewline
GL-3 & 691&71 & 76&61 & 8&22 & 467&01\tabularnewline
GL-4 & 1906&73 & 1048&58 & 9&67 & 992&46\tabularnewline
GL-5 & 3922&59 & 1771&11 & 16&04 & 6495&71\tabularnewline
\hline 
\end{tabular}\caption{Condition numbers of left-preconditioned matrices with preconditioners
Gauss--Seidel, $LD$, and $DU$ with spatial discretization size
$h=2^{-4}$ applied to 2D wave and Klein--Gordon equations with variable
coefficients with $s$-stage Gauss--Legendre. The time-steps are
given by $h_{t}$ = $h^{\frac{1}{s}}$.}
\label{table:left-hyperbolic}
\end{table}

We now repeat the experiments from subsection \ref{subsec:parab-cond},
changing the equations to the diffusion and Klein--Gordon problems
and the timestepper to Gauss--Legendre IRKN. In this case, we see
once again that the variable coefficients have little effect on the
condition numbers. A difference between these problems and the parabolic
problems is that the $GSL$ and $DU$ preconditioners are ineffective
at $s\ge4$ while the $LD$ preconditioner remains effective at all
stage numbers.

\subsection{\label{subsec:GMRES-Performance-of-preconditioned}Performance of
preconditioned GMRES}

The acid test of a preconditioner is of course its performance in
a solve. Now we investigate the effect of our top three preconditioners
in reducing the number of iterations and time required for GMRES to
converge to a solution. As mentioned above, we employ the method of
manufactured solutions in order to examine the relative error norms
rather than just the residuals. For all of the experiments, a relative
residual tolerance for GMRES is set to $10^{-8}$. For each stage-$s$
method (where $s=2\ldots5)$ , we examine the iteration count and
solve time for each preconditioner, in addition to the resulting relative
error. 
\begin{table}
\begin{centering}
\begin{tabular}{cc|cr@{\extracolsep{0pt}.}lr@{\extracolsep{0pt}.}l|cr@{\extracolsep{0pt}.}lc|cr@{\extracolsep{0pt}.}lc}
\hline 
 &  & \multicolumn{5}{c|}{$\mathcal{P}_{GSL}^{-1}\mathcal{A}$} & \multicolumn{4}{c|}{$\mathcal{P}_{LD}^{-1}\mathcal{A}$} & \multicolumn{4}{c}{$\mathcal{P}_{DU}^{-1}\mathcal{A}$}\tabularnewline
\hline 
RIIA & $h$ & it. & \multicolumn{2}{c}{t} & \multicolumn{2}{c|}{err} & it. & \multicolumn{2}{c}{t} & err & it. & \multicolumn{2}{c}{t} & err\tabularnewline
\hline 
\multirow{2}{*}{s=2} & $2^{-5}$ & 11 & 0&21 & 8&0e-10 & 5 & 0&09 & 8.9e-9 & 6 & 0&12 & 1.5e-9\tabularnewline
 & $2^{-6}$ & 8 & 0&83 & 1&1e-8 & 7 & 0&68 & 4.5e-10 & 7 & 0&68 & 2.7e-10\tabularnewline
\hline 
\multirow{2}{*}{s=3} & $2^{-5}$ & 13 & 0&48 & 2&1e-9 & 9 & 0&38 & 6.4e-10 & 10 & 0&41 & 4.7e-9\tabularnewline
 & $2^{-6}$ & 11 & 1&56 & 8&2e-9 & 8 & 1&24 & 1.7e-9 & 10 & 1&47 & 4.9e-9\tabularnewline
\hline 
\multirow{2}{*}{s=4} & $2^{-5}$ & 19 & 0&97 & 1&1e-9 & 10 & 0&57 & 4.3e-9 & 16 & 0&85 & 1.4e-9\tabularnewline
 & $2^{-6}$ & 16 & 2&91 & 4&2e-9 & 10 & 2&05 & 1.5e-9 & 14 & 2&72 & 4.0e-9\tabularnewline
\hline 
\multirow{2}{*}{s=5} & $2^{-5}$ & 19 & 1&26 & 3&9e-9 & 11 & 0&81 & 5.2e-9 & 18 & 1&24 & 5.5e-9\tabularnewline
 & $2^{-6}$ & 17 & 4&07 & 8&0e-9 & 11 & 2&86 & 4.3e-9 & 19 & 4&62 & 1.8e-9\tabularnewline
\hline 
\end{tabular}
\par\end{centering}
\caption{\label{table:GMRES Parablic Beta Zero} Iteration counts, elapsed
time (in seconds), and relative error for left-preconditioned GMRES
with a relative residual tolerance of $10^{-8}$ for a 2D diffusion
equation with variable coefficients with $s$-stage Radau IIA methods
with preconditioners: Gauss--Seidel, LD, and DU. The time-step is
given by $h_{t}=h^{\frac{2}{2s-1}}$. Preconditioners are applied
via a single AMG V-Cycle for all of the block subsolves.}
\end{table}

Rather than apply the preconditioners exactly, we employ the same
strategy as that in \cite{mardal2007,staff2006,masud2020}. Due to
the design of the matrices (either block lower triangular or block
upper triangular) we approximate application of the preconditioners
by solving the block systems using forward or backward substitution.
For each diagonal subsolve, we use a single AMG V-Cycle. This allows
for faster solves with regards to the preconditioner, which translates
to overall faster solve times. For more details on the specifics on
this process, see \cite{masud2020}.

\subsubsection{\label{subsec:GMRES-Diffusion-and-Pennes}Diffusion and Pennes equations}

\begin{table}
\begin{centering}
\begin{tabular}{ccc|cr@{\extracolsep{0pt}.}lc|cr@{\extracolsep{0pt}.}lc|cr@{\extracolsep{0pt}.}lc}
 &  & \multicolumn{1}{c}{} & \multicolumn{4}{c|}{$\mathcal{P}_{GSL}^{-1}\mathcal{A}$} & \multicolumn{4}{c|}{$\mathcal{P}_{LD}^{-1}\mathcal{A}$} & \multicolumn{4}{c}{$\mathcal{P}_{DU}^{-1}\mathcal{A}$}\tabularnewline
\cline{2-15} \cline{3-15} \cline{4-15} \cline{5-15} \cline{7-15} \cline{8-15} \cline{9-15} \cline{11-15} \cline{12-15} \cline{13-15} \cline{15-15} 
 & RIIA & \multicolumn{1}{c}{$h$} & it. & \multicolumn{2}{c}{t} & err. & it. & \multicolumn{2}{c}{t} & err. & it. & \multicolumn{2}{c}{t} & err.\tabularnewline
\cline{2-15} \cline{3-15} \cline{4-15} \cline{5-15} \cline{7-15} \cline{8-15} \cline{9-15} \cline{11-15} \cline{12-15} \cline{13-15} \cline{15-15} 
\multirow{8}{*}{\begin{turn}{90}
Constant Coef.
\end{turn}} & \multirow{2}{*}{s=2} & $2^{-5}$ & 8 & 0&11 & 4.5e-8 & 5 & 0&08 & 3.6e-8 & 6 & 0&09 & 8.7e-9\tabularnewline
 &  & $2^{-6}$ & 8 & 0&55 & 4.1e-8 & 5 & 0&39 & 2.6e-8 & 6 & 0&44 & 4.2e-9\tabularnewline
\cline{2-15} \cline{3-15} \cline{4-15} \cline{5-15} \cline{7-15} \cline{8-15} \cline{9-15} \cline{11-15} \cline{12-15} \cline{13-15} \cline{15-15} 
 & \multirow{2}{*}{s=3} & $2^{-5}$ & 10 & 0&31 & 4.8e-8 & 9 & 0&33 & 8.2e-9 & 8 & 0&26 & 3.8e-9\tabularnewline
 &  & $2^{-6}$ & 14 & 1&53 & 4.5e-9 & 11 & 1&29 & 1.4e-9 & 9 & 0&97 & 1.8e-7\tabularnewline
\cline{2-15} \cline{3-15} \cline{4-15} \cline{5-15} \cline{7-15} \cline{8-15} \cline{9-15} \cline{11-15} \cline{12-15} \cline{13-15} \cline{15-15} 
 & \multirow{2}{*}{s=4} & $2^{-5}$ & 13 & 0&56 & 6.4e-7 & 11 & 0&54 & 2.6e-9 & 15 & 0&70 & 1.3e-8\tabularnewline
 &  & $2^{-6}$ & 19 & 2&86 & 2.9e-9 & 14 & 2&25 & 3.6e-10 & 15 & 2&36 & 5.2e-8\tabularnewline
\cline{2-15} \cline{3-15} \cline{4-15} \cline{5-15} \cline{7-15} \cline{8-15} \cline{9-15} \cline{11-15} \cline{12-15} \cline{13-15} \cline{15-15} 
 & \multirow{2}{*}{s=5} & $2^{-5}$ & 20 & 1&19 & 5.6e-9 & 10 & 0&61 & 5.6e-7 & 18 & 1&09 & 3.4e-8\tabularnewline
 &  & $2^{-6}$ & 23 & 4&47 & 2.9e-9 & 14 & 3&03 & 5.4e-9 & 20 & 4&08 & 1.5e-8\tabularnewline
\hline 
\hline 
\multirow{8}{*}{\begin{turn}{90}
Variable Coef.
\end{turn}} & \multirow{2}{*}{s=2} & $2^{-5}$ & 11 & 0&21 & 7.8e-10 & 5 & 0&09 & 8.9e-9 & 6 & 0&12 & 1.7e-9\tabularnewline
 &  & $2^{-6}$ & 8 & 0&82 & 1.1e-8 & 7 & 0&68 & 4.5e-10 & 7 & 0&67 & 2.8e-10\tabularnewline
\cline{2-15} \cline{3-15} \cline{4-15} \cline{5-15} \cline{7-15} \cline{8-15} \cline{9-15} \cline{11-15} \cline{12-15} \cline{13-15} \cline{15-15} 
 & \multirow{2}{*}{s=3} & $2^{-5}$ & 13 & 0&48 & 2.2e-9 & 9 & 0&38 & 6.4e-10 & 10 & 0&41 & 5.3e-9\tabularnewline
 &  & $2^{-6}$ & 11 & 1&59 & 8.2e-9 & 8 & 1&24 & 1.7e-9 & 10 & 1&48 & 3.9e-9\tabularnewline
\cline{2-15} \cline{3-15} \cline{4-15} \cline{5-15} \cline{7-15} \cline{8-15} \cline{9-15} \cline{11-15} \cline{12-15} \cline{13-15} \cline{15-15} 
 & \multirow{2}{*}{s=4} & $2^{-5}$ & 19 & 0&97 & 1.0e-9 & 10 & 0&57 & 4.3e-9 & 16 & 0&85 & 1.6e-9\tabularnewline
 &  & $2^{-6}$ & 16 & 2&91 & 4.3e-9 & 10 & 2&05 & 1.5e-9 & 15 & 2&85 & 4.7e-9\tabularnewline
\cline{2-15} \cline{3-15} \cline{4-15} \cline{5-15} \cline{7-15} \cline{8-15} \cline{9-15} \cline{11-15} \cline{12-15} \cline{13-15} \cline{15-15} 
 & \multirow{2}{*}{s=5} & $2^{-5}$ & 19 & 1&25 & 4.1e-9 & 11 & 0&81 & 5.2e-9 & 19 & 1&29 & 1.3e-9\tabularnewline
 &  & $2^{-6}$ & 17 & 4&06 & 7.9e-9 & 11 & 2&86 & 4.3e-9 & 20 & 4&81 & 1.4e-9\tabularnewline
\end{tabular}
\par\end{centering}
\caption{Iteration counts, elapsed time (in seconds), and relative error for
left-preconditioned GMRES with a relative residual tolerance of  $10^{-8}$
for a 2D Pennes Equation ($\beta>0)$ with constant coefficients (top)
and variable coefficients (bottom) with $s$-stage Radau IIA methods
with preconditioners: Gauss--Seidel, LD, and DU. The time-step is
given by $h_{t}=h^{\frac{2}{2s-1}}$. }

\label{table: GMRES Parabolic Beta Positive}
\end{table}

For the parabolic problems, we again employ the Radau IIA IRK method
timestepper. Due to the established results for the diffusion problem
in \cite{mardal2007,masud2020,staff2006}, we focus on the Pennes
problem with both constant and variable coefficients. 

Results are shown in Table \ref{table: GMRES Parabolic Beta Positive}.
Iteration counts for all three preconditioners depend weakly on $h$
and on $s$. For fixed $s$, solve time increases by a factor of $\sim3-4$
as $h$ is halved, as expected in a 2D problem with $h$-independent
iteration count. Acceptable errors are achieved, comparable to the
residual tolerance imposed. Once again, the $LD$ preconditioner usually
outperforms the other two. 

\subsubsection{\label{subsec:GMRES-Wave-and-Klein-Gordon}Wave and Klein--Gordon
equations}

\begin{table}
\begin{centering}
\begin{tabular}{ccc|cr@{\extracolsep{0pt}.}lc|cr@{\extracolsep{0pt}.}lc|cr@{\extracolsep{0pt}.}lc}
\hline 
 &  & \multicolumn{1}{c}{} & \multicolumn{4}{c|}{$\mathcal{P}_{GSL}^{-1}\mathcal{A}$} & \multicolumn{4}{c|}{$\mathcal{P}_{LD}^{-1}\mathcal{A}$} & \multicolumn{4}{c}{$\mathcal{P}_{DU}^{-1}\mathcal{A}$}\tabularnewline
\cline{2-15} \cline{3-15} \cline{4-15} \cline{5-15} \cline{7-15} \cline{8-15} \cline{9-15} \cline{11-15} \cline{12-15} \cline{13-15} \cline{15-15} 
 & GL & \multicolumn{1}{c}{$h$} & it. & \multicolumn{2}{c}{t} & err. & it. & \multicolumn{2}{c}{t} & err. & it. & \multicolumn{2}{c}{t} & err.\tabularnewline
\cline{2-15} \cline{3-15} \cline{4-15} \cline{5-15} \cline{7-15} \cline{8-15} \cline{9-15} \cline{11-15} \cline{12-15} \cline{13-15} \cline{15-15} 
\multirow{8}{*}{\begin{turn}{90}
Constant Coef.
\end{turn}} & \multirow{2}{*}{s=2} & $2^{-5}$ & 14 & 0&19 & 5.4e-8 & 7 & 0&10 & 8.5e-8 & 10 & 0&15 & 4.7e-9\tabularnewline
 &  & $2^{-6}$ & 14 & 0&84 & 3.3e-8 & 10 & 0&67 & 6.6e-9 & 10 & 0&66 & 3.6e-9\tabularnewline
\cline{2-15} \cline{3-15} \cline{4-15} \cline{5-15} \cline{7-15} \cline{8-15} \cline{9-15} \cline{11-15} \cline{12-15} \cline{13-15} \cline{15-15} 
 & \multirow{2}{*}{s=3} & $2^{-5}$ & 16 & 0&47 & 2.0e-7 & 8 & 0&29 & 8.9e-9 & 8 & 0&25 & 8.5e-8\tabularnewline
 &  & $2^{-6}$ & 22 & 2&06 & 8.2e-9 & 7 & 0&81 & 1.3e-8 & 8 & 0&81 & 5.1e-8\tabularnewline
\cline{2-15} \cline{3-15} \cline{4-15} \cline{5-15} \cline{7-15} \cline{8-15} \cline{9-15} \cline{11-15} \cline{12-15} \cline{13-15} \cline{15-15} 
 & \multirow{2}{*}{s=4} & $2^{-5}$ & 66 & 2&71 & 3.5e-7 & 29 & 1&24 & 1.1e-8 & 28 & 1&18 & 2.1e-8\tabularnewline
 &  & $2^{-6}$ & 60 & 7&60 & 7.8e-7 & 27 & 3&66 & 8.4e-9 & 20 & 2&60 & 5.6e-7\tabularnewline
\cline{2-15} \cline{3-15} \cline{4-15} \cline{5-15} \cline{7-15} \cline{8-15} \cline{9-15} \cline{11-15} \cline{12-15} \cline{13-15} \cline{15-15} 
 & \multirow{2}{*}{s=5} & $2^{-5}$ & 204 & 11&17 & 2.3e-8 & 31 & 1&80 & 4.5e-8 & 27 & 1&55 & 4.5e-6\tabularnewline
 &  & $2^{-6}$ & 211 & 35&95 & 3.4e-8 & 31 & 5&54 & 3.5e-8 & 28 & 4&95 & 1.5e-6\tabularnewline
\hline 
\hline 
\multirow{8}{*}{\begin{turn}{90}
Variable Coef.
\end{turn}} & \multirow{2}{*}{s=2} & $2^{-5}$ & 12 & 0&15 & 3.8e-8 & 7 & 0&10 & 5.9e-8 & 10 & 0&15 & 6.5e-7\tabularnewline
 &  & $2^{-6}$ & 11 & 0&67 & 2.1e-7 & 6 & 0&46 & 7.5e-8 & 6 & 0&45 & 9.9e-8\tabularnewline
\cline{2-15} \cline{3-15} \cline{4-15} \cline{5-15} \cline{7-15} \cline{8-15} \cline{9-15} \cline{11-15} \cline{12-15} \cline{13-15} \cline{15-15} 
 & \multirow{2}{*}{s=3} & $2^{-5}$ & 19 & 0&71 & 2.0e-8 & 9 & 0&38 & 3.1e-10 & 10 & 0&41 & 5.1e-9\tabularnewline
 &  & $2^{-6}$ & 16 & 2&40 & 2.1e-8 & 6 & 0&95 & 1.4e-8 & 10 & 1&57 & 3.0e-9\tabularnewline
\cline{2-15} \cline{3-15} \cline{4-15} \cline{5-15} \cline{7-15} \cline{8-15} \cline{9-15} \cline{11-15} \cline{12-15} \cline{13-15} \cline{15-15} 
 & \multirow{2}{*}{s=4} & $2^{-5}$ & 61 & 2&53 & 1.6e-7 & 25 & 1&18 & 1.1e-8 & 29 & 1&38 & 9.9e-9\tabularnewline
 &  & $2^{-6}$ & 49 & 7&30 & 9.2e-7 & 20 & 3&63 & 1.2e-8 & 27 & 4&76 & 6.9e-9\tabularnewline
\cline{2-15} \cline{3-15} \cline{4-15} \cline{5-15} \cline{7-15} \cline{8-15} \cline{9-15} \cline{11-15} \cline{12-15} \cline{13-15} \cline{15-15} 
 & \multirow{2}{*}{s=5} & $2^{-5}$ & 164 & 9&05 & 5.8e-8 & 27 & 1&70 & 3.7e-8 & 34 & 2&16 & 1.8e-7\tabularnewline
 &  & $2^{-6}$ & 149 & 29&25 & 5.1e-7 & 35 & 6&11 & 3.4e-8 & 31 & 7&58 & 7.8e-8\tabularnewline
\end{tabular}
\par\end{centering}
\caption{Iteration counts, elapsed time (in seconds), and relative error for
left-preconditioned GMRES with a relative residual tolerance of $10^{-8}$
for a 2D wave equation ($\beta=0)$ with constant coefficients (top)
and variable coefficients (bottom) with $s$-stage Gauss--Legendre
methods with preconditioners: Gauss--Seidel, LD, and DU. The time
step is given by $h_{t}=h^{\frac{1}{s}}$.}

\label{table: GMRES Hyperbolic Beta  0}
\end{table}
\begin{table}
\begin{centering}
\begin{tabular}{ccc|cr@{\extracolsep{0pt}.}lc|cr@{\extracolsep{0pt}.}lc|cr@{\extracolsep{0pt}.}lc}
\hline 
 &  & \multicolumn{1}{c}{} & \multicolumn{4}{c|}{$\mathcal{P}_{GSL}^{-1}\mathcal{A}$} & \multicolumn{4}{c|}{$\mathcal{P}_{LD}^{-1}\mathcal{A}$} & \multicolumn{4}{c}{$\mathcal{P}_{DU}^{-1}\mathcal{A}$}\tabularnewline
\cline{2-15} \cline{3-15} \cline{4-15} \cline{5-15} \cline{7-15} \cline{8-15} \cline{9-15} \cline{11-15} \cline{12-15} \cline{13-15} \cline{15-15} 
 & GL & \multicolumn{1}{c}{$h$} & it. & \multicolumn{2}{c}{t} & err. & it. & \multicolumn{2}{c}{t} & err. & it. & \multicolumn{2}{c}{t} & err.\tabularnewline
\cline{2-15} \cline{3-15} \cline{4-15} \cline{5-15} \cline{7-15} \cline{8-15} \cline{9-15} \cline{11-15} \cline{12-15} \cline{13-15} \cline{15-15} 
\multirow{8}{*}{\begin{turn}{90}
Constant Coef.
\end{turn}} & \multirow{2}{*}{s=2} & $2^{-8}$ & 4 & 0&66 & 1.0e-8 & 3 & 0&64 & 7.8e-9 & 3 & 0&64 & 1.9e-8\tabularnewline
 &  & $2^{-9}$ & 4 & 2&47 & 1.9e-8 & 3 & 2&55 & 1.3e-8 & 3 & 2&58 & 3.2e-8\tabularnewline
\cline{2-15} \cline{3-15} \cline{4-15} \cline{5-15} \cline{7-15} \cline{8-15} \cline{9-15} \cline{11-15} \cline{12-15} \cline{13-15} \cline{15-15} 
 & \multirow{2}{*}{s=3} & $2^{-8}$ & 7 & 2&40 & 4.7e-9 & 5 & 1&97 & 8.6e-9 & 5 & 1&91 & 5.3e-9\tabularnewline
 &  & $2^{-9}$ & 5 & 11&21 & 1.2e-8 & 4 & 11&45 & 6.5e-8 & 4 & 9&68 & 8.1e-9\tabularnewline
\cline{2-15} \cline{3-15} \cline{4-15} \cline{5-15} \cline{7-15} \cline{8-15} \cline{9-15} \cline{11-15} \cline{12-15} \cline{13-15} \cline{15-15} 
 & \multirow{2}{*}{s=4} & $2^{-8}$ & 40 & 14&86 & 4.8e-9 & 11 & 4&75 & 3.1e-8 & 10 & 4&36 & 9.1e-9\tabularnewline
 &  & $2^{-9}$ & 38 & 100&1 & 4.3e-9 & 11 & 34&30 & 6.0e-8 & 15 & 46&87 & 1.7e-9\tabularnewline
\cline{2-15} \cline{3-15} \cline{4-15} \cline{5-15} \cline{7-15} \cline{8-15} \cline{9-15} \cline{11-15} \cline{12-15} \cline{13-15} \cline{15-15} 
 & \multirow{2}{*}{s=5} & $2^{-8}$ & 100 & 47&07 & 1.0e-8 & 14 & 7&45 & 6.5e-8 & 14 & 7&46 & 2.7e-8\tabularnewline
 &  & $2^{-9}$ & 111 & 218&2 & 1.7e-8 & 18 & 39&00 & 1.3e-9 & 16 & 35&44 & 1.3e-8\tabularnewline
\hline 
\hline 
\multirow{8}{*}{\begin{turn}{90}
Variable Coef.
\end{turn}} & \multirow{2}{*}{s=2} & $2^{-8}$ & 5 & 0&74 & 7.9e-9 & 3 & 0&64 & 1.2e-8 & 3 & 0&63 & 2.3e-8\tabularnewline
 &  & $2^{-9}$ & 4 & 2&53 & 5.7e-8 & 3 & 2&69 & 1.8e-8 & 3 & 2&66 & 5.1e-8\tabularnewline
\cline{2-15} \cline{3-15} \cline{4-15} \cline{5-15} \cline{7-15} \cline{8-15} \cline{9-15} \cline{11-15} \cline{12-15} \cline{13-15} \cline{15-15} 
 & \multirow{2}{*}{s=3} & $2^{-8}$ & 12 & 3&82 & 2.7e-9 & 4 & 1&49 & 1.7e-8 & 5 & 1&93 & 7.6e-9\tabularnewline
 &  & $2^{-9}$ & 9 & 11&46 & 3.1e-9 & 4 & 7&22 & 7.2e-8 & 4 & 5&89 & 2.3e-8\tabularnewline
\cline{2-15} \cline{3-15} \cline{4-15} \cline{5-15} \cline{7-15} \cline{8-15} \cline{9-15} \cline{11-15} \cline{12-15} \cline{13-15} \cline{15-15} 
 & \multirow{2}{*}{s=4} & $2^{-8}$ & 40 & 15&13 & 7.2e-9 & 12 & 5&23 & 5.6e-8 & 15 & 6&69 & 1.3e-9\tabularnewline
 &  & $2^{-9}$ & 39 & 95&14 & 4.2e-9 & 12 & 21&10 & 4.7e-8 & 10 & 18&36 & 2.9e-8\tabularnewline
\cline{2-15} \cline{3-15} \cline{4-15} \cline{5-15} \cline{7-15} \cline{8-15} \cline{9-15} \cline{11-15} \cline{12-15} \cline{13-15} \cline{15-15} 
 & \multirow{2}{*}{s=5} & $2^{-8}$ & 100 & 48&32 & 1.8e-8 & 18 & 9&45 & 4.3e-8 & 16 & 8&64 & 2.0e-8\tabularnewline
 &  & $2^{-9}$ & 111 & 221&4 & 1.7e-8 & 17 & 37&08 & 5.8e-8 & 16 & 34&92 & 1.8e-8\tabularnewline
\end{tabular}
\par\end{centering}
\caption{Iteration counts, elapsed time (in seconds), and relative error for
left-preconditioned GMRES to converge to a residual tolerance of $10^{-8}$
for a 2D Klein--Gordon equation ($\beta>0)$ with constant coefficients
(top) and variable coefficients (bottom) with $s$-stage Gauss--Legendre
methods with preconditioners: Gauss--Seidel, LD, and DU. The time-step
is given by $h_{t}=h^{\frac{1}{s}}$.}

\label{table: GMRES Hyperbolic Beta  Positive}
\end{table}

Finally, we investigate the performance of the top three preconditioners
in GMRES as applied to the wave and Klein--Gordon equations with
constant and variable coefficients. As before, the Gauss--Legendre
IRKN method is used for these problem. For the wave equation, we look
at $h=2^{-5}$ and $2^{-6}$, while for the Klein--Gordon equation
we consider $h=2^{-8}$ and $2^{-9}$. Results for the wave equation
are shown in Table \ref{table: GMRES Hyperbolic Beta  0}; those for
the Klein--Gordon equation are shown in \ref{table: GMRES Hyperbolic Beta  Positive}.
For all three preconditioners we see a significant increase in iteration
count and solve time as $s$ is increased, but all show little or
no sensitivity to $h$. Interestingly, unlike in the parabolic problems,
we see the $LD$ and $DU$ preconditioners outperform the $GSL$ preconditioner,
with the $DU$ being slightly better than the $LD$ preconditioner.
Again, tolerable accuracy is achieved.

\section{\label{sec:Conclusions}Conclusions}

In this paper, we developed a unified formulation and analysis of
the stage equations for implicit Runge--Kutta and Runge--Kutta--Nyström
timesteppers applied to a large class of parabolic and hyperbolic
equations of importance in applications. With this unified approach,
we were able to prove the order-optimality of many preconditioners
for these problems, including those in \cite{mardal2007,staff2006}
and \cite{masud2020}. In particular, the order optimality of Rana's
LD preconditioner for any problem is a new result. We also performed
numerical experiments to investigate the dependence on timestep and
mesh size in practice, the effect on these preconditioners on the
spectrum and field of values of the system, and the influence on GMRES
solve time. Since our formulation encompasses problems with variable
coefficients, we also investigated such problems. 

We found in all cases that preconditioner performance was only slightly
influenced by whether the problem's coefficients were variable or
constant. For the parabolic problems with IRK (Radau IIA) timesteppers,
we found results consistent with those of \cite{masud2020}, namely
that the LD preconditioner consistently outperforms the alternatives.
For the hyperbolic problems with IRKN (Gauss--Legendre based) timesteppers
we found that the DU preconditioner --- which hadn't been very effective
on parabolic problems --- was marginally superior to the LD preconditioner,
and both were markedly superior to the Gauss--Seidel preconditioner. 

Finally, we mention that while the class of problems considered here
is large, the restriction to non-negative coefficients $\beta\left(\mathbf{r}\right)\ge0$
may be an issue for linearizations of nonlinear problems, in which
a stage equation will have coefficients involving functions of previous
Newton iterates, which may not respect that restriction. Furthermore,
our analysis does not include the effect of advective terms $\mathbf{b}\cdot\nabla u$,
which will be important in some applications. In \cite{masud2020},
Rana \emph{et al.} found that the LD preconditioner performed well
in the presence of an advective term, but we don't yet have theoretical
understanding of why it works well for those problems. Experimental
and theoretical investigation of these issues will be the subject
of future papers. 

\bibliography{references}

\bibliographystyle{siam}

\newpage{}

\appendix

\section{\label{sec:Proof-of-theorem}Proof of Theorem \ref{thm:iso-A}}

The proof of our Theorem \ref{thm:iso-A} follows the outline of the
proof of the special case in \cite{mardal2007}, with the appropriate
modifications for our more general problem. At each timestep, we solve
the equation $\mathcal{{A}}\mathbf{u}^{n}=\mathbf{f}^{n}$, where
$\mathcal{{A}}$ and $\mathbf{f}^{n}$ were defined in (\ref{eq:cont-A})
and (\ref{eq:RHS}). We introduce the bilinear functional $a\left(\mathbf{u},\mathbf{v}\right):=\left\langle \mathcal{{A}}\mathbf{u},\mathbf{v}\right\rangle _{\mathbf{V}}$.
With few restrictions on the structure of the Butcher coefficient
matrix $A$ in the operator $\mathcal{{A}}$ defined in (\ref{eq:cont-A}),
we cannot assume $a\left(\cdot,\cdot\right)$ is coercive; we therefore
need the Babuska--Aziz theorem \cite{babuska1972,oden2012}. 

Throughout this section, norms and inner products are assumed to be
on $V$ or $\mathbf{V}$ unless otherwise specified. 
\begin{theorem}
\label{thm:(Babuska-Aziz)}(Babuska--Aziz) The linear map $\mathcal{A}:\mathbf{V}\to\mathbf{V}^{*}$
is an isomorphism if the following conditions are satisfied: 
\begin{enumerate}
\item {[}Boundedness{]} There exists a $c_{1}$ independent of $h_{t}$
such that: 
\[
\left|a\left(\mathbf{u},\mathbf{v}\right)\right|\leq c_{1}\left\Vert \mathbf{u}\right\Vert \left\Vert \mathbf{v}\right\Vert ,\quad\forall\mathbf{u},\mathbf{v}\in\mathbf{V}.
\]
\item {[}inf-sup{]} There exists a $c_{2}$ independent of $h_{t}$ such
that 
\[
\sup_{\mathbf{v}\in\mathbf{V}}\frac{\left|a\left(\mathbf{u},\mathbf{v}\right)\right|}{\left\Vert \mathbf{v}\right\Vert }\geq\frac{1}{c_{2}}\left\Vert \mathbf{u}\right\Vert ,\quad\mathbf{u}\in\mathbf{V}.
\]
\item For $\mathbf{v}\in\mathbf{V}\backslash\mathbf{0}$ there exists $\mathbf{u}\in\mathbf{V}$
such that 
\[
a\left(\mathbf{u},\mathbf{v}\right)\neq0.
\]
\end{enumerate}
\end{theorem}

The proof of Theorem \ref{thm:iso-A} therefore reduces to showing
that $\mathcal{{A}}$ satisfies the conditions in Theorem \ref{thm:(Babuska-Aziz)}.
The essential idea is that we can follow the path laid out by \cite{mardal2007},
replacing operations such as $\left\langle \nabla\mathbf{u},\nabla\mathbf{v}\right\rangle $
by $\left\langle \mathbf{u},\mathbf{v}\right\rangle _{\mathbf{H}_{0}^{1}}$;
the logic is the same but the details differ. We require the following
lemma, proved in \cite{nilssen2005}. 
\begin{lemma}
\label{lem:WPD-for-Babuska}Let $A\in\mathbb{{R}}^{s\times s}$ be
a weakly positive definite matrix, and let $C\in\mathbb{{R}}^{s\times s}$
be the associated positive definite matrix of $A$. Then, there exists
$\epsilon>0$ such that for all $\mathbf{x}\in\mathbb{R}^{s}$, we
have 
\begin{align*}
\mathbf{x}^{T}C\mathbf{x} & \geq\epsilon\|\mathbf{x}\|_{2}^{2}\\
\mathbf{x}^{T}CA\mathbf{x} & \geq\epsilon\|\mathbf{x}\|_{2}^{2}.
\end{align*}
\end{lemma}

We now prove Theorem \ref{thm:iso-A}.
\begin{proof}
We begin with the boundedness condition. Let $a_{\text{max}}=\text{max}\left(\max_{i,j}\left|a_{ij}\right|,1\right)$.
Then compute a bound on $\left|a\left(\mathbf{u},\mathbf{v}\right)\right|$
as follows:
\begin{align*}
\left|a\left(\mathbf{u},\mathbf{v}\right)\right| & =\left|\left\langle \mathcal{{A}}\mathbf{u},\mathbf{v}\right\rangle \right|=\left|\left\langle \mathbf{u},\mathbf{v}\right\rangle _{\mathbf{L}^{2}}+h_{t}^{\mu}\sum_{i=1}^{s}\sum_{j=1}^{s}a_{ij}\left\langle u_{j},v_{i}\right\rangle _{H_{0}^{1}}\right|\\
 & \leq\left|\left\langle \mathbf{u},\mathbf{v}\right\rangle _{\mathbf{L}^{2}}\right|+h_{t}^{\mu}\sum_{i,j\le s}\left|a_{ij}\right|\left|\left\langle u_{j},v_{i}\right\rangle _{H_{0}^{1}}\right|\\
 & \leq a_{\text{max}}\left[\left\Vert \mathbf{u}\right\Vert _{\mathbf{L}^{2}}\left\Vert \mathbf{v}\right\Vert _{\mathbf{L}^{2}}+h_{t}^{\mu}\sum_{i,j\le s}\left\Vert u_{j}\right\Vert _{H_{0}^{1}}\left\Vert v_{i}\right\Vert _{H_{0}^{1}}\right]\\
 & \leq a_{\text{max}}s\left[\left\Vert \mathbf{u}\right\Vert _{\mathbf{L}^{2}}\left\Vert \mathbf{v}\right\Vert _{\mathbf{L}^{2}}+h_{t}^{\mu}\left\Vert \mathbf{u}\right\Vert _{\mathbf{H}_{0}^{1}}\left\Vert \mathbf{v}\right\Vert _{\mathbf{H}_{0}^{1}}\right]\\
 & \leq a_{\text{max}}s\left(\left\Vert \mathbf{u}\right\Vert _{\mathbf{L}^{2}}^{2}+h_{t}^{\mu}\left\Vert \mathbf{u}\right\Vert _{\mathbf{H}_{0}^{1}}^{2}\right)^{1/2}\left(\left\Vert \mathbf{v}\right\Vert _{\mathbf{L}^{2}}^{2}+h_{t}^{\mu}\left\Vert \mathbf{v}\right\Vert _{\mathbf{H}_{0}^{1}}^{2}\right)^{1/2}\\
 & \leq a_{\text{max}}s\left\Vert \mathbf{u}\right\Vert _{\mathbf{V}}\left\Vert \mathbf{v}\right\Vert _{\mathbf{V}}.
\end{align*}
Set $c_{1}$ to $a_{\text{max}}s$, and condition 1 (boundedness)
is met. 

We now establish condition 2, the inf-sup condition. If $\mathbf{u}=\mathbf{0}$,
the result is immediate. Otherwise, let $\mathbf{u}\in\mathbf{V}\backslash\mathbf{0}$,
set $\bm{\mathbf{v}}=C^{T}\mathbf{u}$, and define $\epsilon\in\mathbb{R}$
by 
\[
\epsilon:=\min\Big(\min_{\|\mathbf{x}\|=1}\mathbf{x}^{T}C\mathbf{x},\min_{\|\mathbf{x}\|=1}\mathbf{x}^{T}CA\mathbf{x}\Big).
\]
Now compute:
\begin{align*}
\sup_{\mathbf{v}\in\mathbf{V}}\frac{\left\langle \mathcal{{A}}\mathbf{u},\mathbf{v}\right\rangle }{\left\Vert \mathbf{u}\right\Vert \left\Vert \mathbf{v}\right\Vert } & \geq\frac{\left\langle \mathcal{{A}}\mathbf{u},C^{T}\mathbf{u}\right\rangle }{\left\Vert \mathbf{u}\right\Vert \left\Vert C^{T}\mathbf{u}\right\Vert }\\
 & =\frac{\left\langle C\mathcal{A}\mathbf{u},\mathbf{u}\right\rangle }{\left\Vert \mathbf{u}\right\Vert \left\Vert C^{T}\mathbf{u}\right\Vert }\\
 & =\frac{\left\langle C\left(\identity\otimes\mathcal{{I}}+h_{t}^{\mu}A\otimes\mathcal{{K}}\right)\mathbf{u},\mathbf{u}\right\rangle }{\left\Vert \mathbf{u}\right\Vert \left\Vert C^{T}\mathbf{u}\right\Vert }\\
 & =\frac{1}{\left\Vert \mathbf{u}\right\Vert \left\Vert C^{T}\mathbf{u}\right\Vert }\left[\left\langle C\mathbf{u},\mathbf{u}\right\rangle _{\mathbf{L}^{2}}+h_{t}^{\mu}\left\langle CA\mathbf{u},\mathbf{u}\right\rangle _{\mathbf{H}_{0}^{1}}\right].
\end{align*}
From Lemma \ref{lem:WPD-for-Babuska}, we have 
\begin{align*}
\frac{1}{\left\Vert \mathbf{u}\right\Vert \left\Vert C^{T}\mathbf{u}\right\Vert }\left[\left\langle C\mathbf{u},\mathbf{u}\right\rangle _{\mathbf{L}^{2}}+h_{t}^{\mu}\left\langle CA\mathbf{u},\mathbf{u}\right\rangle _{\mathbf{H}_{0}^{1}}\right] & \ge\frac{\epsilon}{\left\Vert \mathbf{u}\right\Vert \left\Vert C^{T}\mathbf{u}\right\Vert }\left[\left\Vert \mathbf{u}\right\Vert _{\mathbf{L}^{2}}^{2}+h_{t}^{\mu}\left\Vert \mathbf{u}\right\Vert _{\mathbf{H}_{0}^{1}}^{2}\right]\\
 & =\frac{\epsilon\left\Vert \mathbf{u}\right\Vert }{\left\Vert C^{T}\mathbf{u}\right\Vert }\ge\frac{\epsilon}{\left\Vert C^{T}\right\Vert }>0.
\end{align*}
Therefore, ${\displaystyle \sup_{\mathbf{v}\in\mathbf{V}}\frac{\left\langle \mathcal{{A}}\mathbf{u},\mathbf{v}\right\rangle }{\left\Vert \mathbf{v}\right\Vert }\geq\frac{1}{c_{2}}\left\Vert \mathbf{u}\right\Vert }$,
where ${\displaystyle c_{2}^{-1}=\epsilon/\left\Vert C^{T}\right\Vert }>0$,
and the inf-sup condition is met. 

Finally, we prove the third condition. As in \cite{mardal2007}: since
$A^{T}$ and $A$ share the same set of eigenvalues, it follows that
$A^{T}$ must also be weakly positive definite (since no eigenvalues
of $A^{T}$ are negative and real). Hence there exists a positive
definite matrix $D$ such that $DA^{T}$ is positive definite as well.
We choose an arbitrary $\mathbf{v}\in\mathbf{V}\backslash\mathbf{0}$
and define $\mathbf{u}=D^{T}\mathbf{v}$. Then 
\begin{align*}
\mathcal{{A}}\left(\mathbf{u},\mathbf{v}\right) & =\left\langle \mathcal{A}D^{T}\mathbf{v},\mathbf{v}\right\rangle =\left\langle \left(\identity\otimes\mathcal{{I}}+h_{t}^{\mu}A\otimes\mathcal{{K}}\right)D^{T}\mathbf{v},\mathbf{v}\right\rangle \\
 & =\left\langle D^{T}\mathbf{v},\mathbf{v}\right\rangle _{\mathbf{L}^{2}}+h_{t}^{\mu}\left\langle AD^{T}\otimes\mathcal{{K}}\mathbf{v},\mathbf{v}\right\rangle _{\mathbf{L}^{2}}\\
 & =\left\langle \mathbf{v},D\mathbf{v}\right\rangle _{L^{2}}+h_{t}^{\mu}\left\langle \mathbf{v},DA^{T}\mathbf{v}\right\rangle _{\mathbf{H}_{0}^{1}}
\end{align*}
which is strictly positive since $D$ and $DA^{T}$ are positive definite,
the operations $\left\langle \cdot,\cdot\right\rangle _{\mathbf{L}^{2}}$
and $\left\langle \cdot,\cdot\right\rangle _{\mathbf{H}_{0}^{1}}$
are positive definite and (at least) positive semidefinite respectively,
and $\mathbf{v\ne0}$. The third condition is met and the proof is
complete.
\end{proof}

\section{\label{sec:Preconditioner-Def}Some Specific Preconditioners}

Here, we define the preconditioners used in Section \ref{sec:NumericalResults}
to precondition the system $\identity\otimes M+h_{t}^{p}A\otimes F$:

\[
\begin{aligned}\mathcal{{P}}_{J}:= & \identity\otimes M+h_{t}^{p}P_{J}\otimes F\\
\mathcal{{P}}_{GSL}:= & \identity\otimes M+h_{t}^{p}P_{GSL}\otimes F\\
\mathcal{{P}}_{TRIU}:= & \identity\otimes M+h_{t}^{p}P_{TRIU}\otimes F,
\end{aligned}
\]

where $P_{J},\,$$P_{GSL,}$and $P_{TRIU}$ are the diagonal, lower
triangular part, and upper triangular part of $A,$ respectively.\\

Let $A=LDU$ be the LDU factorization of $A$. Now set $P_{LD}=LD$
and $P_{DU}=DU$ and define the remaining preconditioners as

\[
\begin{aligned}\mathcal{{P}}_{LD}:= & \identity\otimes M+h_{t}^{p}P_{LD}\otimes F\\
\mathcal{{P}}_{DU}:= & \identity\otimes M+h_{t}^{p}P_{DU}\otimes F.
\end{aligned}
\]

\end{document}